\documentclass[reqno]{amsart}
\usepackage{amssymb,hyperref}

\theoremstyle{plain}
\newtheorem{thm}{Theorem}
\newtheorem{lem}{Lemma}
\newtheorem{cor}{Corollary}

\theoremstyle{definition}
\newtheorem{ex}{Example}

\renewcommand{\Re}{\mathrm{Re}}
\renewcommand{\Im}{\mathrm{Im}}

\title
[Sufficient conditions for strongly Carath\'{e}odory functions]
{Sufficient conditions \\
for strongly Carath\'{e}odory functions}

\author{Hitoshi Shiraishi}
\address{Hitoshi Shiraishi \newline
Department of Mathematics \newline
Kinki University \newline
Higashi-Osaka, Osaka 577-8502, Japan}
\email{shiraishi@math.kindai.ac.jp}

\author{Shigeyoshi Owa}
\address{Shigeyoshi Owa \newline
Department of Mathematics \newline
Kinki University \newline
Higashi-Osaka, Osaka 577-8502, Japan}
\email{owa@math.kindai.ac.jp}

\author{H. M. Srivastava}
\address{H. M. Srivastava \newline
Departmemt of Mathematics and Statistics \newline
University of Victoria \newline
Victoria, British Columbia V8W 3P4, Canada}
\email{harimsri@math.uvic.ca}


\date{}

\begin{document}

\begin{abstract}
Applying the subordination principle for analytic functions in the open unit disk $\mathbb{U}$,
I. H. Kim and N. E. Cho (Comput. Math. Appl. {\bf 59}(2010), 2067-2073) considered some sufficient conditions for Carath\'{e}odory functions.
The purpose of this paper is to discuss some sufficient conditions for the class of strongly Carath\'{e}odory functions in $\mathbb{U}$.
\end{abstract}

\begin{flushleft}
This paper was published in the journal: \\
Panamer. Math. J. {\bf 21} (2011), No. 3, 63--77.
\end{flushleft}
\hrule

\

\

\maketitle

\section{Introduction}

\

Let $\mathcal{H}[a_0,n]$ denote the class of functions $p(z)$ of the form
$$
p(z)
= a_0 + \sum_{k=n}^{\infty} a_k z^k
\qquad (n=1,2,3,\ldots)
$$
which are analytic in the open unit disk $\mathbb{U}=\{z \in \mathbb{C}:|z|<1\}$,
where $a_0\in\mathbb{C}$.

If $p(z)\in\mathcal{H}[a_0,n]$ satisfies
$$
|\arg (p(z)) |
< \frac{\pi}{2}\mu
\qquad (z\in\mathbb{U})
$$
for some real $\mu$ $(0<\mu\leqq1)$,
then we say that $p(z)$ is the strongly Carath\'{e}odory function of order $\mu$ and written by $p(z)\in\mathcal{STP}(\mu)$.

Also,
let $\mathcal{A}_n$ denote the class of functions
$$
f(z)
= z+a_{n+1}z^{n+1}+a_{n+2}z^{n+2}+ \ldots
\qquad(n=1,2,3,\ldots)
$$
that are analytic in $\mathbb{U}$
and $\mathcal{A}\equiv\mathcal{A}_1$.

Further,
let the class $\mathcal{STS}(\mu)$ of $f(z) \in \mathcal{A}_n$ be defined by
$$
\mathcal{STS}(\mu)
= \left\{f(z)\in\mathcal{A}_n:\left|\arg \left(\frac{zf'(z)}{f(z)}\right) \right|<\frac{\pi}{2}\mu,\ 0<\mu\leqq1\right\}
$$
and $\mathcal{S}^*\equiv\mathcal{STS}(1)$. 
This class $\mathcal{STS}(\mu)$ was considered by Shiraishi and Owa \cite{m1ref4} and $f(z) \in \mathcal{STS}(\mu)$ is said to be the strongly starlike function of order $\mu$ in $\mathbb{U}$.

Let $f(z)$ and $g(z)$ be analytic in $\mathbb{U}$.
Then $f(z)$ is said to be subordinate to $g(z)$ if there exists an analytic function $w(z)$ in $\mathbb{U}$ satisfying $w(0)=0$, $|w(z)| < 1$ $(z\in\mathbb{U})$ and such that $f(z)=g(w(z))$.
We denote this subordination by
$$
f(z)
\prec g(z)
\qquad (z\in\mathbb{U}).
$$
In particular, if $g(z)$ is univalent in $\mathbb{U}$, then the subordination
$$
f(z)
\prec g(z)
\qquad (z\in\mathbb{U})
$$
is equivalent to $f(0)=g(0)$ and $f(\mathbb{U}) \subset g(\mathbb{U})$ (cf. \cite{d7ref2}).

Denote by $\mathcal{Q}$ the class of functions $q(z)$ that are analytic and injective on $\overline{\mathbb{U}} \setminus \mathbb{E}(q)$,
where
$$
\mathbb{E}(q)
=\left\{ \zeta \in \partial\mathbb{U} : \lim_{z\rightarrow\zeta} \{ q(z) \} = \infty \right\},
$$
and such that $q'(\zeta) \neq 0$ $(\zeta\in\partial\mathbb{U}\setminus\mathbb{E}(q))$.

Futher,
let the subclass of $\mathcal{Q}$ for which $q(0)=a_0$ be denoted by $\mathcal{Q}(a_0)$.

\

To discuss our problems,
we need the following lemma due to Miller and Mocanu \cite{d7ref2}.

\

\begin{lem} \label{d7lem1} \quad
Let $q(z)\in\mathcal{Q}(a_0)$ and let $h(z)\in\mathcal{H}[a_0,n]$ with $h(z) \not\equiv a_0$.
If $h(z) \nprec q(z)$,
then there exist points $z_0\in\mathbb{U}$ and $\zeta_0\in\partial\mathbb{U}\setminus\mathbb{E}(q)$ for which
$$
h(z_0)=q(\zeta_0)
$$
and
$$
z_0h'(z_0)=m\zeta_0q'(\zeta_0)
\qquad (m \geqq n \geqq 1).
$$
\end{lem}

\section{Conditions for the strongly  Carath\'{e}odory function}

\

Applying Lemma \ref{d7lem1},
we derive

\

\begin{thm} \label{d7thm1} \quad
Let $g(z)$ be analytic in $\mathbb{U}$ with
\begin{equation} \label{d7thm1eq01}
A=\inf_{z\in\mathbb{U}} \{ \Re(g(z))\cos\alpha - |\Im(g(z))\sin\alpha|  \}>0
\qquad (z\in\mathbb{U})
\end{equation}
for some $-\dfrac{\pi}{2} < \alpha < \dfrac{\pi}{2}$.
If $p(z)\in\mathcal{H}[1,n]$ satisfies $p(z) \not\equiv 1$ and
$$
\Re(p(z)+g(z)zp'(z))
>\frac{1}{2nA}((\cos\alpha+2nA)\sin^2\alpha-n^2A^2\cos\alpha)
\qquad (z\in\mathbb{U}),
$$
then
$$
|\arg (p(z))|
<\frac{\pi}{2}-|\alpha|
\qquad (z\in\mathbb{U}).
$$
\end{thm}

\

\begin{proof} \quad
First,
let us define the function $h_1(z)$ by
\begin{equation} \label{d7thm1eq02}
h_1(z)
= e^{i\alpha} p(z)
\qquad (z\in\mathbb{U})
\end{equation}
and the function $q_1(z)$ by
\begin{equation} \label{d7thm1eq03}
q_1(z)
= \frac{e^{i\alpha}+\overline{e^{i\alpha}}z}{1-z}
\qquad (z\in\mathbb{U})
\end{equation}
for $-\dfrac{\pi}{2}<\alpha<\dfrac{\pi}{2}$.

Then,
we see that $h_1(z)$ and $q_1(z)$ are analytic in $\mathbb{U}$ with
$$
h_1(0)
=q_1(0)
=e^{i\alpha}
\in\mathbb{C}
$$
and
$$
q_1(\mathbb{U})
=\{ w\in\mathbb{C}:\Re(w)>0 \}.
$$

Now we suppose that $h_1(z)$ is not subordinate to $q_1(z)$.
Then, Lemma \ref{d7lem1} shows us that,
there exist points $z_1\in\mathbb{U}$ and $\zeta_1\in\partial\mathbb{U}\setminus \{ 1 \}$ such that
\begin{equation} \label{d7thm1eq04}
h_1(z_1)
=q_1(\zeta_1)
=i\rho_1
\qquad (\rho_1\in\mathbb{R})
\end{equation}
and
\begin{equation} \label{d7thm1eq05}
z_1h_1'(z_1)
=m\zeta_1 q_1'(\zeta_1)
\qquad (m\geqq n \geqq 1).
\end{equation}

Here we note that
\begin{equation} \label{d7thm1eq06}
\zeta_1
=q_1^{-1}(h_1(z_1))
=\frac{h_1(z_1)-e^{i\alpha}}{h_1(z_1)+\overline{e^{i\alpha}}}
\end{equation}
and
\begin{equation} \label{d7thm1eq07}
\zeta_1q_1'(\zeta_1)
=-\frac{\rho_1^2-2\rho_1\sin\alpha+1}{2\cos\alpha}
\equiv \sigma_1(\rho_1)
<0.
\end{equation}

For such $z_1\in\mathbb{U}$ and $\zeta_1\in\partial\mathbb{U}\setminus \{ 1 \}$, we obtain
\begin{align*}
&\Re(p(z_1)+g(z_1)z_1p'(z_1)) \\
&= \Re(e^{-i\alpha}h_1(z_1)+g(z_1)e^{-i\alpha}z_1h_1'(z_1)) \\
&= \Re(e^{-i\alpha}q_1(\zeta_1)+g(z_1)e^{-i\alpha}m\zeta_1q_1'(\zeta_1)) \\
&= \Re( e^{-i\alpha}i\rho_1+g(z_1)e^{-i\alpha}m\sigma_1(\rho_1)) \\
&= \rho_1\sin\alpha + m(\Re(g(z_1))\cos\alpha-\Im(g(z_1))\sin\alpha)\sigma_1(\rho_1) \\
&\leqq \rho_1\sin\alpha + nA\sigma_1(\rho_1) \\
&= \rho_1\sin\alpha - nA\frac{\rho_1^2-2\rho_1\sin\alpha+1}{2\cos\alpha} \\
&= -\frac{nA}{2\cos\alpha} \left( \rho_1-\frac{\sin\alpha(\cos\alpha+nA)}{nA} \right)^2 + \frac{\sin^2\alpha(\cos\alpha+2nA)-n^2A^2\cos\alpha}{2nA} \\
&\leqq \frac{1}{2nA} (\sin^2\alpha(\cos\alpha+2nA)-n^2A^2\cos\alpha),
\end{align*}
where $A$ is given by (\ref{d7thm1eq01}).
This evidently contradicts the assumption of Theorem \ref{d7thm1}.
Therefore we obtain
\begin{equation} \label{d7thm1eq08}
\Re(h_1(z)) = \Re(e^{i\alpha}p(z))>0
\qquad (z\in\mathbb{U})
\end{equation}
for $-\dfrac{\pi}{2}<\alpha<\dfrac{\pi}{2}$.

Next,
let us put
\begin{equation} \label{d7thm1eq09}
h_2(z)
= e^{-i\alpha} p(z)
\qquad (z\in\mathbb{U})
\end{equation}
and
\begin{equation} \label{d7thm1eq10}
q_2(z)
= \frac{e^{-i\alpha}+\overline{e^{-i\alpha}}z}{1-z}
\qquad (z\in\mathbb{U})
\end{equation}
for $-\dfrac{\pi}{2}<\alpha<\dfrac{\pi}{2}$.

Then,
we see that the functions $h_2(z)$ and $q_2(z)$ are analytic in $\mathbb{U}$ with
$$
h_2(0)
=q_2(0)
=e^{-i\alpha}
\in\mathbb{C}
$$
and
$$
q_2(\mathbb{U})
=\{ w\in\mathbb{C}:\Re(w)>0 \}
=q_1(\mathbb{U}).
$$

If we suppose that $h_2(z)$ is not subordinate to $q_2(z)$,
then Lemma \ref{d7lem1} gives us that,
there exist points $z_2\in\mathbb{U}$ and $\zeta_2\in\partial\mathbb{U}\setminus \{ 1 \}$ such that
\begin{equation} \label{d7thm1eq11}
h_2(z_2)
=q_2(\zeta_2)
=i\rho_2
\qquad (\rho_2\in\mathbb{R})
\end{equation}
and
\begin{equation} \label{d7thm1eq12}
z_2h_2'(z_2)
=m\zeta_2 q_2'(\zeta_2)
\qquad (m\geqq n \geqq 1).
\end{equation}

Futher, we note that
\begin{equation} \label{d7thm1eq13}
\zeta_2
=q_2^{-1}(h_1(z_2))
=\frac{h_2(z_2)-e^{-i\alpha}}{h_2(z_2)+\overline{e^{-i\alpha}}}
\end{equation}
and
\begin{equation} \label{d7thm1eq14}
\zeta_2q_2'(\zeta_2)
=-\frac{\rho_2^2+2\rho_2\sin\alpha+1}{2\cos\alpha}
\equiv \sigma_2(\rho_2)
<0.
\end{equation}

For such $z_2\in\mathbb{U}$ and $\zeta_2\in\partial\mathbb{U}\setminus \{ 1 \}$ , we see that
\begin{align*}
&\Re(p(z_2)+g(z_2)z_2p'(z_2)) \\
&= \Re(e^{i\alpha}h_2(z_2)+g(z_2)e^{i\alpha}z_2h_2'(z_2)) \\
&= \Re(e^{i\alpha}q_2(\zeta_2)+g(z_2)e^{i\alpha}m\zeta_2q_2'(\zeta_2)) \\
&= \Re( e^{i\alpha}i\rho_2+g(z_2)e^{i\alpha}m\sigma_2(\rho_2) ) \\
&= -\rho_2\sin\alpha + m(\Re(g(z_2))\cos\alpha+\Im(g(z_2))\sin\alpha)\sigma_2(\rho_2) \\
&\leqq -\rho_2\sin\alpha + nA\sigma_2(\rho_2) \\&= -\rho_2\sin\alpha - nA\frac{\rho_2^2+2\rho_2\sin\alpha+1}{2\cos\alpha} \\
&= -\frac{nA}{2\cos\alpha} \left( \rho_2+\frac{\sin\alpha(\cos\alpha+nA)}{nA} \right)^2 + \frac{\sin^2\alpha(\cos\alpha+2nA)-n^2A^2\cos\alpha}{2nA} \\
&\leqq \frac{1}{2nA} (\sin^2\alpha(\cos\alpha+2nA)-n^2A^2\cos\alpha),
\end{align*}
where $A$ is given by (\ref{d7thm1eq01}).
This also contradicts the assumption of Theorem \ref{d7thm1}.
Therefore we have that
\begin{equation} \label{d7thm1eq15}
\Re(h_2(z)) = \Re(e^{-i\alpha}p(z))>0
\qquad (z\in\mathbb{U})
\end{equation}
for $-\dfrac{\pi}{2}<\alpha<\dfrac{\pi}{2}$.

Hence, combining inequalities (\ref{d7thm1eq08}) and (\ref{d7thm1eq15}),
we complete the proof of Theorem \ref{d7thm1}.
\end{proof}

\

Taking $\dfrac{\pi}{2}\mu=\dfrac{\pi}{2}-\alpha$ $(0<\mu\leqq1)$ in Theorem \ref{d7thm1},
we obtain 

\

\begin{cor} \label{d7cor1} \quad
Let $g(z)$ be analytic in $\mathbb{U}$ with
$$
A=\inf_{z\in\mathbb{U}} \left\{ \Re(g(z)) \sin\frac{\pi}{2}\mu - \left| \Im(g(z))\cos\frac{\pi}{2}\mu \right| \right\}>0
\qquad (z\in\mathbb{U})
$$
for some $0 < \mu \leqq 1$.
If $p(z)\in\mathcal{H}[1,n]$ satisfies $p(z) \not\equiv 1$ and
$$
\Re(p(z)+g(z)zp'(z))
>\frac{1}{2nA} \left( \left( \sin\frac{\pi}{2}\mu+2nA \right) \cos^2\frac{\pi}{2}\mu-n^2A^2\sin\frac{\pi}{2}\mu \right)
\quad (z\in\mathbb{U}),
$$
then $p(z)\in\mathcal{STP}(\mu)$
\end{cor}

\

Putting
$$
p(z)=\dfrac{zf'(z)}{f(z)}=1+na_{n+1}z^n+\ldots
\qquad (z\in\mathbb{U})
$$
for $f(z)\in\mathcal{A}_n$ in Corollary \ref{d7cor1},
we have

\

\begin{cor} \label{d7cor2} \quad
Let $g(z)$ be analytic in $\mathbb{U}$ with
$$
A=\inf_{z\in\mathbb{U}} \left\{ \Re(g(z)) \sin\frac{\pi}{2}\mu - \left| \Im(g(z))\cos\frac{\pi}{2}\mu \right| \right\}>0
\qquad (z\in\mathbb{U})
$$
for some $0 < \mu \leqq 1$.
If $f(z)\in\mathcal{A}_n$ satisfies $\dfrac{zf'(z)}{f(z)} \not\equiv 1$ and
\begin{align*}
&\Re \left( \frac{zf'(z)}{f(z)} +g(z)\frac{zf'(z)}{f(z)} \left( 1 -\frac{zf'(z)}{f(z)} +\frac{zf''(z)}{f'(z)} \right) \right) \\
&>\frac{1}{2nA} \left( \left( \sin\frac{\pi}{2}\mu+2nA \right) \cos^2\frac{\pi}{2}\mu-n^2A^2\sin\frac{\pi}{2}\mu \right)
\qquad (z\in\mathbb{U}),
\end{align*}
then
$f(z)\in\mathcal{STS}(\mu)$.
\end{cor}

\

Considering $n=1$ and $0\leqq\alpha<\dfrac{\pi}{2}$,
we get 

\

\begin{cor} \label{d7cor3} \quad
Let $g(z)$ be analytic in $\mathbb{U}$ with
$$
A=\inf_{z\in\mathbb{U}} \{ \Re(g(z))\cos\alpha - \Im(g(z))\sin\alpha  \}>0
\qquad (z\in\mathbb{U})
$$
for some $0 \leqq \alpha < \dfrac{\pi}{2}$.
If $p(z)\in\mathcal{H}[1,1]$ satisfies $p(z) \not\equiv 1$ and
$$
\Re(p(z)+g(z)zp'(z))
>\frac{1}{2A}((\cos\alpha+2A)\sin^2\alpha-A^2\cos\alpha)
\qquad (z\in\mathbb{U}),
$$
then
$$
|\arg (p(z)) |
<\frac{\pi}{2}-\alpha
\qquad (z\in\mathbb{U}).
$$
\end{cor}

\

We try to show an example of Theorem \ref{d7thm1}.

\

\begin{ex} \label{d7ex1} \quad
Let us consider the function
$$
p(z)
= 1+kz^n
\qquad (n=1,2,3,\ldots)
$$
for $k \leqq \dfrac{n^2+6n-9}{4n(4n+3)}$.
Then, it is easy to observe that
$p(z)$ is analytic in $\mathbb{U}$ and maps $\mathbb{U}$ onto the disk with the center at $p(0)=1$ and the radius $k$.

Thus,
we see that
$$
|\arg (p(z))|
<\sin^{-1}\left( \frac{n^2+6n-9}{4n(4n+3)} \right)
<\frac{\pi}{4}
\qquad (z\in\mathbb{U}).
$$

Further,
if we consider $g(z)=1+\dfrac{1}{3}z$ and $\alpha=\dfrac{\pi}{4}$,
we obtain that
\begin{align*}
A
&=\inf_{z\in\mathbb{U}} \{ \Re(g(z))\cos\alpha - |\Im(g(z))\sin\alpha|  \} \\
&=\frac{1}{3\sqrt[]{2}}
>0
\qquad (z\in\mathbb{U}).
\end{align*}

This gives us that
$$
\frac{1}{2nA}((\cos\alpha+2nA)\sin^2\alpha-n^2A^2\cos\alpha)
= \frac{-n^2+6n+9}{12n}.
$$

On the other hand,
we also have
\begin{align*}
\Re(p(z)+g(z)zp'(z))
&= \Re \left( 1+(1+n)kz^n+\frac{n}{3}kz^{n+1} \right) \\
&> 1-(1+n)k-\frac{n}{3}k \\
&\geqq \frac{-n^2+6n+9}{12n}
\qquad (z\in\mathbb{U}).
\end{align*}

Thus,
in view of Theorem \ref{d7thm1}, we conclude that
$$
|\arg (p(z))|
< \frac{\pi}{4}
\qquad (z\in\mathbb{U}).
$$
\end{ex}

\

We also derive

\

\begin{thm} \label{d7thm2} \quad
If $p(z)\in\mathcal{H}[1,n]$ satisfies $p(z)\neq0$ $(z \in \mathbb{U})$, $p(z)\not\equiv1$ and
$$
-\frac{\sqrt[]{2n\cos^2\alpha+n^2}+n\sin\alpha}{\cos\alpha}
< \Im \left( p(z)+\frac{zp'(z)}{p(z)} \right)
< \frac{\sqrt[]{2n\cos^2\alpha+n^2}-n\sin\alpha}{\cos\alpha}
$$
for all $z\in\mathbb{U}$ and some $0\leqq\alpha<\dfrac{\pi}{2}$,
then
$$
|\arg (p(z))|
< \frac{\pi}{2}-\alpha
\qquad (z\in\mathbb{U}).
$$
\end{thm}

\

\begin{proof} \quad
We define the functions $h_1(z)$ by  (\ref{d7thm1eq02}) and $q_1(z)$ by (\ref{d7thm1eq03}) for $-\dfrac{\pi}{2}<\alpha<\dfrac{\pi}{2}$.

If $h_1(z)$ is not subordinate to $q_1(z)$,
then there exist points $z_1\in\mathbb{U}$ and $\zeta_1\in\partial\mathbb{U}\setminus \{ 1 \}$ satisfying (\ref{d7thm1eq04}) and (\ref{d7thm1eq05}).

Using equations from (\ref{d7thm1eq02}) to (\ref{d7thm1eq07}),
we have
\begin{align*}
\Im \left( p(z_1)+\frac{z_1p'(z_1)}{p(z_1)} \right)
&= \Im \left( e^{-i\alpha}h_1(z_1)+\frac{z_1h_1'(z_1)}{h_1(z_1)} \right) \\
&= \Im \left( e^{-i\alpha}q_1(\zeta_1)+\frac{m\zeta_1q_1'(\zeta_1)}{q_1(\zeta_1)} \right) \\
&= \rho_1\cos\alpha - \frac{m\sigma_1(\rho_1)}{\rho_1}
\qquad (\rho_1\in\mathbb{R}\setminus \{ 0 \})
\end{align*}
for such $z_1\in\mathbb{U}$ and $\zeta_1\in\partial\mathbb{U}\setminus \{ 1 \}$.

For the case $\rho_1>0$,
since $\sigma_1(\rho_1)<0$ and $m \geqq n$,
we obtain that
\begin{align*}
\rho_1\cos\alpha - \frac{m\sigma_1(\rho_1)}{\rho_1}
&\geqq \rho_1\cos\alpha - \frac{n\sigma_1(\rho_1)}{\rho_1} \\
&= \frac{\rho_1^2(2\cos^2\alpha+n)-2n\rho_1\sin\alpha+n}{2\rho_1\cos\alpha}
\qquad (\rho_1\in\mathbb{R}\setminus \{ 0 \}).
\end{align*}

Since the function
$$
g_1(\rho_1)
=\frac{\rho_1^2(2\cos^2\alpha+n)-2n\rho_1\sin\alpha+n}{2\rho_1\cos\alpha}
\qquad (\rho_1\in\mathbb{R}\setminus \{ 0 \})
$$
takes the minimum value at $\rho_1^*$ given by
$$
\rho_1^*
= \sqrt[]{\frac{n}{2\cos^2\alpha+n}},
$$
we have
\begin{align*}
\rho_1\cos\alpha - \frac{m\sigma_1(\rho_1)}{\rho_1}
&\geqq g_1(\rho_1^*) \\
&= \frac{\sqrt[]{2n\cos^2\alpha+n^2}-n\sin\alpha}{\cos\alpha}
\qquad (\rho_1\in\mathbb{R}\setminus \{ 0 \}),
\end{align*}
whch is the contradiction for the assumption of Theorem \ref{d7thm2}.

For the case $\rho_1<0$,
we put $\rho_1=-\rho_1'$ $(\rho_1'>0)$.

Then, using the same method mentioned above,
we have
\begin{align*}
\rho_1\cos\alpha - \frac{m\sigma_1(\rho_1)}{\rho_1}
&\leqq -\rho_1'\cos\alpha + \frac{n\sigma_1(-\rho_1')}{\rho_1'} \\
&= -\frac{{\rho_1'}^2(2\cos^2\alpha+n)+2n\rho_1'\sin\alpha+n}{2\rho_1'\cos\alpha} \\
&\equiv g_1(-\,\rho_1') \\
&\leqq g_1 \left(-\, \sqrt[]{\frac{n}{2\cos^2\alpha+n}} \right) \\
&= -\frac{\sqrt[]{2n\cos^2\alpha+n^2}+n\sin\alpha}{\cos\alpha},
\end{align*}
which contradicts the assumption of Theorem \ref{d7thm2}.
Hence, we have
\begin{equation} \label{d7thm2eq1}
\Re(e^{i\alpha}p(z))>0
\qquad (z\in\mathbb{U})
\end{equation}
for $0\leqq\alpha<\dfrac{\pi}{2}$.

Next,
considering the functions $h_2(z)$ defined by (\ref{d7thm1eq09}) and $q_2(z)$ defined by (\ref{d7thm1eq10}) and using the similar method as the above,
we also get
\begin{equation} \label{d7thm2eq2}
\Re(e^{-i\alpha}p(z))>0
\qquad (z\in\mathbb{U})
\end{equation}
for $0\leqq\alpha<\dfrac{\pi}{2}$.

Therefore, making use of (\ref{d7thm2eq1}) and (\ref{d7thm2eq2}),
we have complete the proof of Theorem \ref{d7thm2}.
\end{proof}

\

Considering $\dfrac{\pi}{2}\mu=\dfrac{\pi}{2}-\alpha$ $(0<\mu\leqq1)$ in Theorem \ref{d7thm2},
we obtain Corollary \ref{d7cor4}.

\

\begin{cor} \label{d7cor4} \quad
If $p(z)\in\mathcal{H}[1,n]$ satisfies $p(z)\neq0$ $(z \in \mathbb{U})$, $p(z)\not\equiv1$ and
$$
-\frac{\sqrt[]{2n\sin^2\dfrac{\pi}{2}\mu+n^2}+n\cos\dfrac{\pi}{2}\mu}{\sin\dfrac{\pi}{2}\mu}
< \Im \left( p(z)+\frac{zp'(z)}{p(z)} \right)
< \frac{\sqrt[]{2n\sin^2\dfrac{\pi}{2}\mu+n^2}-n\cos\dfrac{\pi}{2}\mu}{\sin\dfrac{\pi}{2}\mu}
$$
for all $z\in\mathbb{U}$ and some $0<\mu\leqq1$,
then $p(z)\in\mathcal{STP}(\mu)$
\end{cor}

\

Putting
$$
p(z)=\dfrac{zf'(z)}{f(z)}=1+na_{n+1}z^n+\ldots
\qquad (z\in\mathbb{U})
$$
for $f(z)\in\mathcal{A}_n$ in Corollary \ref{d7cor4},
we have

\

\begin{cor} \label{d7cor5} \quad
If $f(z)\in\mathcal{A}_n$ satisfies $\dfrac{zf'(z)}{f(z)}\neq0$ $(z \in \mathbb{U})$, $f(z)\not\equiv1$ and
$$
-\frac{\sqrt[]{2n\sin^2\dfrac{\pi}{2}\mu+n^2}+n\cos\dfrac{\pi}{2}\mu}{\sin\dfrac{\pi}{2}\mu}
< \Im \left( 1+\frac{zf''(z)}{f'(z)} \right)
< \frac{\sqrt[]{2n\sin^2\dfrac{\pi}{2}\mu+n^2}-n\cos\dfrac{\pi}{2}\mu}{\sin\dfrac{\pi}{2}\mu}
$$
for all $z\in\mathbb{U}$ and some $0<\mu\leqq1$,
then $f(z)\in\mathcal{STS}(\mu)$
\end{cor}

\

Taking $n=1$ in Theorem \ref{d7thm2},
we obtain the following corollary due to Kim and Cho \cite{d7ref1}.

\

\begin{cor} \label{d7cor6} \quad
If $p(z)\in\mathcal{H}[1,1]$ satisfies $p(z)\neq0$ $(z \in \mathbb{U})$, $p(z)\not\equiv1$ and
$$
-\frac{\sqrt[]{2\cos^2\alpha+1}+\sin\alpha}{\cos\alpha}
< \Im \left( p(z)+\frac{zp'(z)}{p(z)} \right)
< \frac{\sqrt[]{2\cos^2\alpha+1}-\sin\alpha}{\cos\alpha}
\qquad(z \in \mathbb{U})
$$
for some $0\leqq\alpha<\dfrac{\pi}{2}$,
then
$$
|\arg (p(z))|
< \frac{\pi}{2}-\alpha
\qquad (z\in\mathbb{U}).
$$
\end{cor}

\

Finaly,
we derive

\

\begin{thm} \label{d7thm3} \quad
If $p(z)\in\mathcal{H}[1,n]$ satisfies $p(z)\neq0$ $(z \in \mathbb{U})$, $p(z)\not\equiv1$ and
$$
\left| p(z)+\frac{zp'(z)}{p(z)}-1 \right|
< \left( \frac{n}{2}+1 \right) |p(z)|\cos\alpha
\qquad (z\in\mathbb{U})
$$
for some $-\dfrac{\pi}{2}<\alpha<\dfrac{\pi}{2}$,
then
$$
|\arg (p(z)) |
< \frac{\pi}{2}-|\alpha|
\qquad (z\in\mathbb{U}).
$$
\end{thm}

\

\begin{proof} \quad
Let us define the function $h_1(z)$ by
$$
h_1(z)
= \frac{e^{i\alpha}}{p(z)}
\qquad (z\in\mathbb{U})
$$
and $q_1(z)$ be the function defined as in the proof of Theorem \ref{d7thm1}.
If $h_1(z)$ is not subordinate to $q_1(z)$,
then there exist points $z_1\in\mathbb{U}$ and $\zeta_1\in\partial\mathbb{U}\setminus \{ 1 \}$ satisfying (\ref{d7thm1eq04}) and (\ref{d7thm1eq05}).

Using equations from (\ref{d7thm1eq03}) to (\ref{d7thm1eq07}),
we have
\begin{align*}
&\frac{\left| p(z_1)+\dfrac{z_1p'(z_1)}{p(z_1)}-1 \right|}{|p(z_1)|} \\
&= |1-e^{-i\alpha}z_1h_1'(z_1)-e^{-i\alpha}h_1(z_1)| \\
&= |h_1(z_1)+z_1h_1'(z_1)-e^{-i\alpha}| \\
&= |q_1(\zeta_1)+m\zeta_1q_1'(\zeta_1)-e^{-i\alpha}| \\
&= |(m\sigma_1(\rho_1)-\cos\alpha)+i(\rho_1-\sin\alpha)| \\
&= ((m\sigma_1(\rho_1)-\cos\alpha)^2+(\rho_1-\sin\alpha)^2)^{\frac{1}{2}} \\
&= \left( \left( m\frac{\rho_1^2-2\rho_1\sin\alpha+1}{2\cos\alpha}+\cos\alpha \right)^2 + (\rho_1-\sin\alpha)^2 \right)^{\frac{1}{2}} \\
&\geqq  \left(  \left( n\frac{\rho_1^2-2\rho_1\sin\alpha+1}{2\cos\alpha}+\cos\alpha \right)^2 + (\rho_1-\sin\alpha)^2 \right)^{\frac{1}{2}} \\
&= \left( \left( n\frac{(\rho_1-\sin\alpha)^2+\cos^2\alpha}{2\cos\alpha}+\cos\alpha \right)^2 + (\rho_1-\sin\alpha)^2 \right)^{\frac{1}{2}} \\
&\geqq \left( \frac{n}{2}+1 \right) \cos\alpha,
\end{align*}
which is the contradiction for the assumption of Theorem \ref{d7thm3}.
Hence we have
\begin{equation} \label{d7thm3eq1}
\Re(h_1(z))=\Re \left( \frac{e^{i\alpha}}{p(z)} \right)
> 0
\qquad (z\in\mathbb{U}).
\end{equation}

Next,
we consider the function $h_2(z)$ defined by
$$
h_2(z)
= \frac{e^{-i\alpha}}{p(z)}
\qquad (z\in\mathbb{U}).
$$

Then, by virtue of (\ref{d7thm1eq10}), we also obtain that
\begin{equation} \label{d7thm3eq2}
\Re(h_2(z))=\Re \left( \frac{e^{-i\alpha}}{p(z)} \right)
> 0
\qquad (z\in\mathbb{U}).
\end{equation}

Therefore, in view of (\ref{d7thm3eq1}) and (\ref{d7thm3eq2}),
we prove Theorem \ref{d7thm3}.
\end{proof}

\

Taking $\dfrac{\pi}{2}\mu=\dfrac{\pi}{2}-\alpha$ $(0<\mu\leqq1)$ in Theorem \ref{d7thm3},
we obtain

\

\begin{cor} \label{d7cor7} \quad
If $p(z)\in\mathcal{H}[1,n]$ satisfies $p(z)\neq0$ $(z \in \mathbb{U})$, $p(z)\not\equiv1$ and
$$
\left| p(z)+\frac{zp'(z)}{p(z)}-1 \right|
< \left( \frac{n}{2}+1 \right) |p(z)|\sin\frac{\pi}{2}\mu
\qquad (z\in\mathbb{U})
$$
for some $0<\mu\leqq1$,
then $p(z)\in\mathcal{STP}(\mu)$.
\end{cor}

\

Considering the function
$$
p(z)=\dfrac{zf'(z)}{f(z)}=1+na_{n+1}z^n+\ldots
\qquad (z\in\mathbb{U})
$$
for $f(z)\in\mathcal{A}_n$ in Corollary \ref{d7cor7},
we have the following corollary.

\

\begin{cor} \label{d7cor8} \quad
If $f(z)\in\mathcal{A}_n$ satisfies $\dfrac{zf'(z)}{f(z)}\neq0$ $(z \in \mathbb{U})$, $f(z)\not\equiv1$ and
$$
\left| \frac{zf''(z)}{f'(z)} \right|
< \left( \frac{n}{2}+1 \right) \left| \frac{zf'(z)}{f(z)} \right| \sin\frac{\pi}{2}\mu
\qquad (z\in\mathbb{U})
$$
for some $0<\mu\leqq1$,
then $f(z)\in\mathcal{STS}(\mu)$.
\end{cor}

\

Putting $n=1$ and $0\leqq\alpha<\dfrac{\pi}{2}$ in Theorem \ref{d7thm3},
we get Corollary \ref{d7cor9} due to Kim and Cho \cite{d7ref1}.

\

\begin{cor} \label{d7cor9} \quad
If $p(z)\in\mathcal{H}[1,1]$ satisfies $p(z)\neq0$ $(z \in \mathbb{U})$, $p(z)\not\equiv1$ and
$$
\left| p(z)+\frac{zp'(z)}{p(z)}-1 \right|
< \frac{3}{2}|p(z)|\cos\alpha
\qquad (z\in\mathbb{U})
$$
for some $0\leqq\alpha<\dfrac{\pi}{2}$,
then
$$
|\arg (p(z))|
< \frac{\pi}{2}-\alpha
\qquad (z\in\mathbb{U}).
$$
\end{cor}

\

\end{document}